\documentclass[12pt,amsmath]{article}
\usepackage{amssymb}
\usepackage{amscd}
\newtheorem{thm}{Theorem}[section]
\newtheorem{lem}[thm]{Lemma}

\newtheorem{rem}[thm]{Remark}
\newtheorem{cor}[thm]{Corollary}

\newtheorem{defn}[thm]{Definition}

\def\proof{{\sc Proof}.\ }
\def\endproof{\hfill$\square$}
\def\Mat{\mathop{\mathrm {Mat}}\nolimits}
\def\M{\mathop{\mathrm {M}}\nolimits}
\def\Hom{\mathop{\mathrm {Hom}}\nolimits}
\def\HP{\mathop{\mathrm {HP}}\nolimits}
\def\C{$\check {\rm C}$}
\def\Cech{$\check {\rm C}$ech\ }
\def\Tot{\mathop{\mathrm {Tot}}\nolimits}

\def\Spec{\mathop{\mathrm {Spec}}\nolimits}

\em
\begin{document}

\title{$\mathbf Z_2$-graded $\check{\rm C}$ech Cohomology in Noncommutative Geometry\footnote{The work was supported in part by Vietnam National Project of Research in  Fundamental Sciences and The Abdus Salam ICTP, UNESCO.}}
\author{Do Ngoc Diep}
\maketitle

\begin{abstract}
The $\mathbf Z_2$-graded $\check{\rm C}$ech cohomology theory is considered in the framework of noncommutative geometry over complex number field and in particular the homotopy invariance and Morita invariance are proven. In some special case we deduce an isomorphism between this noncommutative theory and the classical $\mathbf Z_2$-graded $\check{\rm C}$ech cohomology theory. 
\end{abstract}

\section{Introduction}
Let us fix in this paper the field of complex numbers as the ground field for algebras, modules, etc. 
 
The general $\check{\rm C}$ech cohomology presheaf $\check{\mathcal H}^q$ for an arbitrary presheaf $X$ on a category $\mathcal C$, $\check{H}^q(X,G) = \check{\mathcal H}^q(G)(X)$ was introduced in \cite{sga4.2}.

\par
The idea of $\check{\rm C}$ech cohomology for noncommutative geometry was appeared in \cite{rosenberg_kontsevich}, \cite{rosenberg}. In this paper we use this idea to define the corresponding (periodic) $\mathbf Z_2$-graded \Cech theory.

\par
We prove the homotopy invariance and Morita invariance of  $\check{\rm C}$ech cohomology in the framework of noncommutative geometry.
Our main result is based on a detailed analysis of the structure of C*-algebras.
A crucial observation is the fact that for C*-algebras the category of *-representations defines  exactly the C*-algebra itself, by the well-known Gelfand-Naimark-Segal Theorem. From this we can deduce the Morita invariance and homotopy invariance of the \C ech cohomology, the same properties of periodic cyclic homology of the C*-algebra. Since our result is valid not only for C*-algebras, we work in the general context of a noncommutative algebra over complex numbers.

\par The paper is organized as follows.
Taking the \Cech cohomology in place of the de Rham theory in the periodic cyclic  theory of A. Connes, in \S2 we define the $\mathbf Z_2$-graded \Cech cohomology theory. Then in section 3 we prove two important properties of the theory  as the homotopy invariance and Morita invariance. In the last section 4 we deduce also some kind of Connes-Hochschild-Kostant-Rosenberg theorem. This lets us see a clear relation with the classical case of commutative algebras and ordinary \Cech cohomology theory.
\par {\it Notations}: We follow  the notations from \cite{sga4.2} and \cite{orlov}

\section{$\mathbf Z_2$-graded $\check{\bf C}$ech cohomology} 
\subsection{Grothendieck topos}
\par
The main purpose of this section is to formulate and define the functor of (periodic cyclic) $\mathbf Z_2$-graded \C ech cohomology. The well-known periodic cyclic homology is based on the cyclic homology theory of A. Connes, which is an algebraic framework of the $\mathbf Z_2$-graded de Rham cohomology theories. In the algebraic context, it was defined by J. Cuntz and D. Quillen in terms of $\mathcal X$-complexes and it has become a new chapter of noncommutative algebraic geometry. The most general Gronthendieck algebraic geometry is purely based in terms of categories. In the generic case this turns out to the algebraic version of the \C ech cohomology in place of de Rham cohomology theories. We follows the work of Orlov \cite{orlov} in particular to formulate the theory. The main references are \cite{sga4.2} and \cite{orlov}.
\par
Many of our results could be obtained in the fields of other characteristics, but we restrict ourselves to the complex case.
\par
Let us denote by $\mathcal C$ a fixed category and $\mathbf {Set}$ the category of sets. Any contravariant functor $X$ from $\mathcal C$ to $\mathbf {Set}$ is called a {\sl presheaf of sets} and the category of all presheaves of sets on  $\mathcal C$ is denoted by $\hat\mathcal C$. The category $\mathcal C$ can be considered as a subcategory of $\hat\mathcal C$, consisting of representable functors $h = h_R : \mathcal C \to \mathbf {Set}.$
If $R$ is an object of $\mathcal C$, then there is a natural isomorphism $\Hom_{\hat \mathcal C}(h_R,X) = X(R).$ For any object $X\in \hat\mathcal C$ the {\sl category over $X$} is the category of pairs $(R, \Phi)$, where $R$ is an object of $\mathcal C$ and $\Phi\in X(R)$, and is denoted by $\mathcal C/X$. 
\par
Recall that a {\it sieve} in the category $\mathcal C$ is a full subcategory $\mathcal D \subseteq \mathcal C$ such that
any object of $\mathcal C$ for which there exists a morphism from it to some object in $\mathcal D$ is contained in $Obj(\mathcal D)$.
A sieve on $R$ is nothing more than a {\sl subpresheaf} of $R$ in the category $\hat\mathcal C$
\par
A {\it Grothendieck topology} $\mathcal T$ on a category $\mathcal C$ is defined by giving for each object $R$ in $\mathcal C$ a set $J(R)$ of the so called {\it covering sieves} satisfying the following axioms:
\begin{enumerate}
\item[(T1)] For any object $R$ the maximal sieve $\mathcal C/R$ is in $J(R)$.
\item[(T2)] If $T \in J(R)$ and $f: S \to R$ a morphism in $\mathcal C$, then the induced sieve $$\CD f^*(T) := \{ U @>\alpha>> S | f\alpha \in T \}\endCD$$ is in $J(S)$.
\item[(T3)] If $T\in J(R)$ is a covering sieve and $U$ is a sieve on $R$ such that $f^*(U) \in J(S)$ for all $f: S\to  R$ in $T$, then $U\in J(R).$
\end{enumerate}
A {\it Grothendieck site} $\Phi = (\mathcal C,\mathcal T)$ is a category $\mathcal C$ and equipped with a Grothendieck topology $\mathcal T$.
\par
It is reasonable to remind the Jacobson topology on the set of all representations of a group  or Zariski topology on algebraic varieties.
\par
For the categories with fiber product, a Grothendieck topology can be given by a {\it Grothendieck pretopology} which is defined by giving for each object $R$ in $\mathcal C$ a family $Cov(R)$ of morphism to $R$ such that
\begin{enumerate}
\item[(P1)] For any family $\{ R_\alpha \to R\}_{\alpha\in I}$ in $Cov(R)$ and $S\to R$ a morphism of $\mathcal C$, the fiber product family $R_\alpha \times_R S \to S$ is also in $Cov(R)$.
\item[(P2)] If $\{R_\alpha \to R\}_{\alpha\in R}$ is in $Cov(R)$ and $\{R_{\beta_\alpha} \to R_\alpha\}_{\beta_\alpha \in J_\alpha}$ is in $Cov(R_\alpha)$ for each $\alpha\in I$, then the total family $\{R_\gamma\to R\}_{\gamma\in \coprod_{\alpha\in I} J_\alpha}$ is in $Cov(R)$.  
\item[(P3)] The trivial family $\{ id_R : R \to R\}$ is in $Cov(R)$.
\end{enumerate}
Any Grothendieck pretopology $P$ on $\mathcal C$ generates a Grothendieck topology $\mathcal T$ such that a sieve is covering in $\mathcal T$ if and only if it contains some covering family in $P$.
\par
The topos on the category of functors from a category $\mathcal C$ to another one $\mathcal D$ is defined by the usual rule.

\subsection{The standard cosimplicial complex of a continuous functor}

We define in this subsection the $\mathbf Z_2$-graded \C ech cohomology theory.
Let us recall the definition of the \C ech cohomology with coefficients in a sheaf $M$.
Let $\mathcal U = (U_i \to X)_{i\in I}$ be a covering sieve of $X$ in the category $\mathcal C$. Suppose that the cover has the property that all fiber product and pushout diagrams exist, see (\cite{sga4.2}, Exp. 4). Denote $\mathcal A_\bullet$ the associated standard simplicial complex:
{\scriptsize
$$\leqno{\mathcal A_\bullet :} \quad
\dots \begin{array}{c} \longrightarrow \\ \longrightarrow \\ \longrightarrow\\ \longrightarrow \end{array} \coprod_{(i,j,k)\in I^3} A(U_i\times U_j \times U_k)\begin{array}{c} \longrightarrow \\ \longrightarrow \\ \longrightarrow \end{array} \coprod_{(i,j)\in I^2} A(U_i\times U_j) \begin{array}{c} \longrightarrow \\ \longrightarrow \end{array}\coprod_{i\in I} A(U_i)$$}

Let us consider again a ringed cite (category) $(\mathcal C, A)$, $\hat\mathcal C$ the topos of presheaves on $\mathcal C$, $\varepsilon : \mathcal C \hookrightarrow \hat\mathcal C$ the canonical functor associating to each object $X\in \mathcal C$ the functor $h_X$, presented by $X$. Define $H^q(h_X,M) = H^q(X,M)$, see (\cite{sga4.2}, Exp IV, 2.3.1), as derived functor of the projective limit functor : $$H^q(S,M) := R^q\lim_{\leftarrow\atop C/S} M|_S$$

For any $A$-module $M$, denote $C^\bullet(\mathcal U,M) := \Hom_A(\mathcal A_\bullet,M)$ :
$$\leqno{C^\bullet(\mathcal U,M):\quad }
\prod_{i\in I} M(U_i) \begin{array}{c} \rightarrow\\ \rightarrow\end{array}
\prod_{(i,j)\in I^2} M(U_i\times U_j) \begin{array}{c} \longrightarrow\\
\longrightarrow\\ \longrightarrow\end{array} \dots $$

One defines $H^q(\mathcal U,M) := H^q(C^\bullet(\mathcal U,M))$. If $R$ is a covering sieve generalized by the family $\mathcal U$ then
$$H^q(\mathcal U,M) \cong H^q(R,M)$$ and the functor $H^q(\mathcal U,.)$ commutes with restriction of scalars, see (\cite{sga4.2}, Exp. IV, Proposition 2.3.4).

One defines, (\cite{sga4.2}, Exp. IV, 2.4)
$$\mathcal H^0(M)(X) := H^0(X,M) = M(X),$$
$$\mathcal H^q(M)(X) := H^q(M,M).$$

For an arbitrary presheaf $G$ of $A$-module $G$, the groups $H^q(R,G)$ are called the {\it \Cech cohomology with respect to the covering sieve $R$, with coefficients in $G$.}

For a sheaf $M$ of $A$-modules over $\mathcal C$, the group 
$$H^q(\mathcal U, M):= H^q(\mathcal U, \mathcal H^0(M))$$ is defined as the {\it \Cech cohomology group of the sheaf $M$ with respect to he the cover $\mathcal U$.}

One has also
$$\check\mathcal H^0(G)(X) = \lim_{\longrightarrow\atop{R\hookrightarrow X}} G(R)$$ and therefore
$$H^q(G)(X) = \lim_{\longrightarrow\atop\ {R\hookrightarrow X}} H^q(R,G)$$ which is called {\it the presheaf of \Cech cohomology}

Define $$\check H^q(X,G) := \check\mathcal H^q(G)(X),$$ one has also
$$\check H^q(X,G) = \lim_{\longrightarrow\atop \mathcal U} H^q(\mathcal U, G).$$

For a sheaf $M$ of $A$-modules, one has
$$\check H^q(X,M) = \check H^q(X, \mathcal H^0(M)),\quad \check \mathcal H^q(M) = \check\mathcal H^q(\check\mathcal H^0(M)).$$
The groups $\check H^q(X,M)$ are called the {\it \Cech cohomology groups of the sheaf $M$.}

\subsection{The periodic cyclic bicomplex}
\begin{lem}[The action of $\mathbf Z_{k+1}$]
There is a natural action of the cyclic group $\mathbf Z_{k+1}$ on the \C ech cohomology cochain complex $C^\bullet(\mathcal U,M)$ associated with a covering $\mathcal U$.
\end{lem}
\proof The action of the cyclic group $\mathbf Z_{k+1}$ is defined a cyclic permutation of indices of $U's$, i.e.
$$(\lambda M)(U_{i_0} \times \dots \times U_{i_k}) := M(U_{i_k}\times U_{i_0}\times\dots \times U_{i_{k-1}}).$$

It is not hard also to see that for a covering sieve $\mathcal U = \{ U_\alpha \to X\}$ there is a natural isomorphism
$$M(U_{i_0} \times \dots \times U_{i_k}) \cong M(U_{i_0}) \otimes \dots \otimes M(U_{i_k}).$$ 
Therefore the \C ech cohomology complex becomes the cyclic complex for $$\lim_{\rightarrow\atop \mathcal U} M(\mathcal U).$$
\endproof

\begin{cor}[Hochschild differentials and Cyclic operations]
The well-known Hochschild differentials $b'$ and $b$ and Connes cyclic operators $\lambda$, $N= 1 + \lambda + \dots + \lambda^k$, $s$ are well-defined on $\mathbf Z_2$-graded \Cech cocycles
\end{cor}

\begin{defn}[Periodic  bicomplex]{\rm
Let $(\mathcal C,A)$ be a ringed $U$-cite, $\hat\mathcal C$ the topos of sheaves, $M$ a sheaf of $A$-modules. Then the bicomplex}

{\tiny
$$\CD
@. \vdots @. \vdots @.  \vdots\\
@. @A-b'AA  @AbAA @A-b'AA @.\\
\dots @>1-\lambda>> \prod_{(i_1,i_2,i_3)\in I^3} M(U_{i_1}\times U_{i_2} \times U_{i_3})@>N>> \prod_{(i_1,i_2,i_3)\in I^3} M(U_{i_1}\times U_{i_2} \times U_{i_3})@>1-\lambda>> \dots\\
@. @A-b'AA  @AbAA @A-b'AA @.\\
\dots @>1-\lambda>> \prod_{(i,j)\in I^2}M(U_i\times U_j) @>N>> \prod_{(i,j)\in I^2}M(U_i\times U_j) @>1-\lambda>>\dots\\
@. @A-b'AA  @AbAA @A-b'AA @.\\
\dots @>1-\lambda>> \prod_{i\in I}M(U_i) @>N>> \prod_{i\in I}M(U_i) @>1-\lambda>> \dots\\
\endCD$$}
{\rm is well-defined and is called the {\it (periodic) \Cech bicomplex}.}
\end{defn}

\begin{defn}[The total complex and $\mathbf Z_2$-graded \Cech cohomology]{\rm 
The {\it associated total complex} of which is defined as
$$\Tot C(\mathcal U,M)^\pm := \prod_{i+j =\pm (mod \; 2)} C^{i,j},$$ where $C^{i,j} := \prod_{(i_1,\dots,i_k)\in I^k} M(U_{i_1}\times\dots \times U_{i_k})$ and $\pm =$ $ev$ (even) or $od$ (odd).

The cohomology of this total complex is called {\it the $\mathbf Z_2$-graded \C ech cohomology of $M$} and denoted by $\mathbf Z_2\check H^*(\mathcal U,M)$ and $\mathbf Z_2\check H(X,M) := \lim_{\longrightarrow\atop \mathcal U} \mathbf Z_2\check H(\mathcal U,M)$. It can be also realized as the cohomology of the total complex related with the process of passing through direct limits, i.e. the direct limit bi-complex
{\tiny
$$\CD
@. \vdots @. \vdots @. \vdots  \\
@. @A-b'AA  @AbAA @A-b'AA @.\\
\dots @>1-\lambda>> \prod_{(i_1,i_2,i_3)\in I^3}\lim_{\longrightarrow\atop \mathcal U} M(U_{i_1}\times U_{i_2} \times U_{i_3})@>N>> \prod_{(i_1,i_2,i_3)\in I^3} \lim_{\longrightarrow\atop \mathcal U}M(U_{i_1}\times U_{i_2} \times U_{i_3})@>1-\lambda>> \dots\\
@. @A-b'AA  @AbAA @A-b'AA \\
\dots @>1-\lambda>> \lim_{\longrightarrow\atop \mathcal U}\prod_{(i,j)\in I^2}M(U_i\times U_j) @>N>> \lim_{\longrightarrow\atop \mathcal U}\prod_{(i,j)\in I^2}M(U_i\times U_j) @>1-\lambda>>  \dots\\
@. @A-b'AA  @AbAA @A-b'AA \\
\dots @>1-\lambda>> \lim_{\longrightarrow\atop \mathcal U}\prod_{i\in I}M(U_i) @>N>> \lim_{\longrightarrow\atop \mathcal U}\prod_{i\in I}M(U_i) @>1-\lambda>> \dots\\
\endCD$$}
For a C*-algebra $A$, we define it $\mathbf Z_2$-graded \Cech cohomology $\mathbf Z_2\check H(A) = \mathbf Z_2\check H(\mathcal A)$ as the $\mathbf Z_2$-graded \Cech cohomology of the category of *-representations of $A$.
}\end{defn}

\begin{rem}
In the first periodic bi-complex without direct limits, all the horizontal lines are acylic, but it is in general not the case for the second periodic bi-complex with direct limits. 
\end{rem}

\section{Homotopy invariance and Morita invariance}
We prove in this section two main properties of the (periodic cyclic) $\mathbf Z_2$-graded \Cech cohomology theory: homotopy invariance and Morita invariance, which make the theory easier to compute and being a generalized homology theory.





\begin{defn}[Chain Homotopy of functors]
Let us consider two functors $F,G: \mathcal C \to \mathcal D$. Denote the corresponding chain functors between complexes by $\{F_n\}, \{G_n\}$, where $F_n,G_n : C_n \to D_n$ for complexes
$$\CD
0 @>>> C_0 @>>> C_1 @>>> C_2 @>>> \dots\\
@. @VF_0,G_0VV @VF_1, G_1VV @VF_2,G_2 VV\\
0 @>>> D_0 @>>> D_1 @>>> D_2 @>>> \dots \endCD$$
We say that $F$ and $G$ are {\it chain homotopic} if there exist  augmentation functors $s_n: C_n \to D_{n-1}$ such that for all $n$
$$F_n - G_n = s_n\circ \partial_{n-1} + \partial_n\circ s_{n+1}.$$
\end{defn}

\begin{lem}
Two functors $F,G: \mathcal A \to \mathcal B$ are homotopic if and only if for any covering sieve $\mathcal U = (U_i \to X)_{i\in I}$, there exists a chain homotopy of chain complexes
$$\coprod_{i\in I} A(U_i) \begin{array}{c}\longleftarrow\\ \longleftarrow\end{array} \coprod_{i,j\in I\times I} A(U_i \times_X U_j)\begin{array}{c}\longleftarrow\\ \longleftarrow\\ \longleftarrow\end{array} \dots$$
and 
$$\coprod_{i\in I} B(U_i) \begin{array}{c}\longleftarrow\\ \longleftarrow\end{array} \coprod_{i,j\in I\times I} B(U_i \times_X U_j)\begin{array}{c}\longleftarrow\\ \longleftarrow\\ \longleftarrow\end{array} \dots$$
\end{lem}

\begin{rem}
In the case of smooth manifolds the chain complex homotopy is realized by integration of the so called Cartan homotopy formula for the Lie derivative
$$L_\xi = \imath(\xi)\circ d + d \circ \imath(\xi)$$
between de Rham complexes. 
\end{rem}


\begin{lem} Let $A$ and $B$ be C*-algebras, and 
let $\mathcal A$ (resp. $\mathcal B$) be the category of $*$-modules . Then the categories $\mathcal A$ and $\mathcal B$ are homotopic one-to-another if and only if the two algebras $A$ and $B$ are homotopic.
\end{lem}
\proof Because of the Gelfand-Naimark-Segal theorem, the C*-algebras are exactly defined by the category of *-representations, the category of *-representations of $B \otimes C[0,1]$ is isomorphic to the category of *-representations of $B$.

\endproof

\begin{lem} Let $A$ be a C*-algebra and $\mathcal A$ the category of *-representations (i.e. $A$-modules) of $A$, then
$$\mathbf Z_2\check{H}^*(\mathcal A) \cong \HP^*(A).$$ 
\end{lem}
\proof
Let us consider affine covering sieve $U_i \to X=\hat{A}=\Spec A$, the dual object of $A$.
\endproof

\begin{thm}[Homotopy Invariance]
Let $\varphi_t : A \to B , t\in I=[0,1]$ be a homotopy of algebras, then 
$$\mathbf Z_2\check{H}^*(A) \cong \mathbf Z_2\check{H}^*(B).$$
\end{thm}
\par
\proof Let $\mathcal A$ (resp. $\mathcal B$) be the category of $A$-modules (resp., $B$-modules).

{\it Step 1.} Change the homotopy by a piecewise-linear homotopy in the space of functors from the category $\mathcal A$ to the category $\mathcal B$.

\par
{\it Step 2.} A piecewise-linear homotopy gives rise to a chain complex homotopy.
\par
{\it Step 3.}
Two chain complex homotopical functors induces the same isomorphism of \C ech cohomology groups. It is an easy consequence from the results of homological algebra: If $F_n$ and $G_n$ are chain complex homotopic than the induced morphisms satisfies
$$F^*_n - G^*_n = \partial_{n-1}^* \circ s^*_n + s^*_{n+1} \circ \partial_n^*.$$
The second summand is a zero morphism on cohomology and the first summand is a boundary. The sum on the right is therefore a zero morphism. 
\endproof

\begin{lem}
Let us denote by $\Mat_n(\mathbb C)$ the algebra of all square $n\times n$- matrices with complex entries. Then we have a natural isomorphism
$$\mathbf Z_2\check H^*(\Mat_n(\mathbb C)) \cong \mathbf Z_2\check H^*(\mathbb C).$$
\end{lem}
\proof
Every complex matrix can be homotopic to a unitary one. Then, every unitary matrix can be by conjugation reduced to a diagonal matrix of complex numbers of module 1. Every elementary block (in this case, diagonal element)$[e^{i\theta}]$ is  homotopic to identity $[1]$ by the classical homotopy $[e^{i\theta t}]_{0\leq t \leq 1}$, i.e.
$$\left[\begin{array}{cc} \cos \theta & -\sin \theta\\ \sin \theta & \cos \theta\end{array}\right] \sim I = \left[\begin{array}{cc} 1 & 0 \\ 0 & 1\end{array}\right]. $$
\endproof

\begin{lem}[Adjoint functors]\label{der}
There is a natural equivalence of functors $\Hom$ and $\otimes$:
$$\Hom(R\otimes \M_n(\mathbb C),M) \cong \Hom(R,\Hom(\M_n(\mathbb C),M)).$$
\end{lem}
\proof
This isomorphism of functors is a particular case of the general adjointness between $\Hom$ and $\otimes$ in homological algebra.
\endproof

\begin{lem}\label{dim}
There is a natural isomorphism of derived functors
$$R^q\Hom(\M_n(\mathbb C),M) \cong R^q\Hom(\mathbb C,M).$$
\end{lem}
\proof
It is an easy exercise from homological algebra. \endproof

\begin{thm}[Morita Invariance]
$$\mathbf Z_2\check{H}^*(A \otimes \Mat_n(\mathbf C)) \cong \mathbf Z_2\check{H}^*(A).$$
\end{thm}
\proof
Let us remark that $A \to A \otimes \Mat_n(\mathbf C)$ is a fiber bundle. Now apply the Grothendieck's Leray-Serre spectral sequence for this fibration.
Following the previous lemmas \ref{der} and \ref{dim}, there is a natural isomorphism of functors
$$R^p\Hom(R\otimes \M_n(\mathbb C), R^q(M) \cong R^p\Hom(R,R^q\Hom(\M_n(\mathbb C),M)),$$ which are the $E^2$ term of a  Leray-Serre spectral sequence converging to the \C ech cohomology.
\endproof

\begin{cor}
The $\mathbf Z_2$-graded \C ech cohomology theory is a generalized cohomology theory.   
\end{cor}

\section{Comparison with the classical $\check{\bf C}$ech cohomology theory}

In this section  we show that a generalization of the Connes-Hochschild-Kostant-Rosenberg Theorem can be obtained easily.

\begin{thm}
Let $A$ be a stable continuous C*-algebra with spectrum a smooth compact manifold $X$, in fact $A=C(X,\mathcal K(P))$ is the algebra of continuous sections of a smooth, locally trivial bundle $\mathcal K(P) := P \times_{PU} \mathcal K$ on $X$ with fibre the algebra $\mathcal K$ of compact operators on a separable Hilbert space associated to a principal $PU$ bundle $P$ on $X$ via the adjoint action of $PU$ on $\mathcal K$.  
Let $\delta(P)\in H^3(X;\mathbf Z)$ be the Dixmier-Douady invariant, that classifies such algebras $A$ and $c(P)$ some closed $3$-form on $X$, that presents the class ${2\pi i}\delta(P)$ in the real cohomology. 
Let $\mathcal A$ be the category of all *-representations of the C*-algebra $A$.Then the $\mathbf Z_2$-graded \Cech cohomology $\mathbf Z_2\check H^*(\mathcal A)$ is isomorphic to the de Rham cohomology $H^*(X;c(P))$ which is isomorphic to the classical $\mathbf Z_2$-graded \Cech cohomology $\mathbf Z_2\check H^*(X; c(P))$. 
\end{thm}
\proof 
In this situation, the $\mathbf Z_2$-graded \Cech cohomology of the category $\mathcal A$ is isomorphic to the Connes periodic cyclic homology $HP_*(C^\infty(X,\mathcal L^1(P)))$, where $C^\infty(X,\mathcal L^1(P))$ is consisting of all smooth section of the sub-bundle $\mathcal L^1(P) = P \times_{PU} \mathcal L^1$ of $\mathcal K(P)$ with fibre the algebra $\mathcal L^1$ of trace class operators on the Hilbert space with the same structure groups $PU$,
see \cite{mathai_stevenson} for a more detailed proof in the language of periodic cyclic homology. \endproof 

\section*{Acknowledgments}
The main part of this work was done while the author was visiting The Abdus Salam ICTP in Trieste, Italy. The author is grateful to ICTP and in particular would like to express his sincere thanks to Professor Dr. Le Dung Trang for invitation and support.

{\parindent=0pt \sc Institute of Mathematics, Vietnam Academy of Science and Technology, 18 Hoang Quoc Viet Road, Cau Giay District, 10307 Hanoi, Vietnam}\\ {\tt Email: dndiep@math.ac.vn}
\end{document}